normalized tight frames and of symmetric orthogonalization of bases by
Hilbert-Schmidt property 
\documentclass[12pt]{amsart}
\theoremstyle{plain}
\newtheorem{theorem}{Theorem}[section]
\newtheorem{corollary}[theorem]{Corollary}
\newtheorem{example}[theorem]{Example}

\newtheorem{proposition}[theorem]{Proposition}
\newtheorem{lemma}[theorem]{Lemma}

\theoremstyle{definition}
\newtheorem{definition}[theorem]{Definition}

\theoremstyle{remark}

\begin{document}

\title{Symmetric approximation of frames and bases in Hilbert spaces}
\author[M.~Frank]{Michael Frank}
\address{Dept.~Mathematics, Univ.~of Houston, Houston, TX 77204-3476, U.S.A.}
\email{frank@math.uh.edu, frank@mathematik.uni-leipzig.de}
\thanks{The first and second authors were supported in part by an NSF grant.}
\author[V.~I.~Paulsen]{Vern I.~Paulsen}
\address{Dept.~Mathematics, Univ.~of Houston, Houston, TX 77204-3476, U.S.A.}
\email{vern@math.uh.edu}
\author[T.~R.~Tiballi]{Terry R.~Tiballi}
\address{Dept.~Mathematics, SUNY at Oswego, Oswego, NY 13126, U.S.A.}
\email{tiballi@oswego.edu}
\begin{abstract}
We consider existence and uniqueness of symmetric approximation of frames by
normalized tight frames and of symmetric orthogonalization of bases by
orthonormal bases in Hilbert spaces.  A crucial role is played by the
Hilbert-Schmidt property of a certain operator related to the initial frame
or basis.
\end{abstract}
\maketitle

Given a Hilbert space $H$ and a linearly independent set of vectors $\{f\}_1^n$
in $H$ the Gram-Schmidt process is traditionally used as a means of creating
an orthonormal set of vectors $\{ \mu_i \}_1^n$ from $\{ f_i \}_1^n$ in such
a way that ${\rm span}\{ f_i \}_1^k = {\rm span}\{ \mu_i \}_1^k$ for all
$1 \leq k \leq n$. This process is inherently order-dependent in that a
reordering of $\{ f_i \}_1^n$ generally results in an entirely new orthonormal
set $\{ \mu_i \}_1^n$. For some problems in quantum chemistry, it is desirable
to treat the vectors $\{ f_1,f_2, ... , f_n \}$ simultaneously in a method that
is order-independent instead of successively as in the Gram-Schmidt process.
It is also desirable in such applications to orthonormalize the set $\{ f_i
\}_1^n$ in such a way that the sum $\sum_{j=1}^n \|\mu_j-f_j\|^2$ is minimal.
The resulting set $\{ \mu_i \}_1^n$ is called the {\it symmetric} or {\it
L\"owdin orthogonalization of} $\{ f_i \}_1^n$.

The original work on this subject was
done by Per-Olov L\"owdin, a quantum chemist, in the late 1940's (\cite{Loew}).
In later publications by Jerome A.~Goldstein and by the coauthors J.~G.~Aiken,
J.~A.~Erdos and M.~Levy (\cite{AEG1,AEG2,GoLe}) this process was extensively
studied. Also, symmetric orthogonalization has been linked to the construction
of optimal algorithms for finding matrix inverse square roots of positive
invertible matrices and for computing the principal square roots of invertible
normal matrices in computer science (\cite{Phi,She1,She2,Lak}).

\smallskip 
In his Ph.D.~thesis the third author investigated the existence and
uniqueness of symmetric orthogonalization of countably infinite linearly
independent sets $\{ f_i \}_1^\infty$ in Hilbert spaces $H$. Introducing the
operator $F: l_2 \to H$, $F(e_i)=f_i$ for any $i \in \mathbb N$
and for the standard orthonormal basis $\{ e_i \}_1^\infty$ of $l_2$, he
showed the existence of symmetric orthogonalizations of $\{ f_i \}_1^\infty$
if $(I-|F|)$ is a Hilbert-Schmidt operator on $l_2$ and the dimension of the
kernel of $F$ is less than or equal to the dimension of the orthogonal
complement of the range of $F$. Uniqueness can be guaranteed if and only if
the dimension of the kernel of $F$ is zero, \cite[Th.~4]{Tib}. An alternative
characterization of a symmetric orthogonalization $\{ \nu_i \}_1^\infty$ of
the set $\{ f_i \}_1^\infty$ is the finiteness of the sum $\sum_{j=1}^\infty
\|\nu_j-f_j\|^2$ and its equality to the square of the Hilbert-Schmidt norm
of the operator $(I-|F|)$, \cite[Th.~2, 3 and Cor.~5]{Tib}. However, the
identity of the Hilbert (sub-)spaces generated by $\{ \nu_i \}_1^\infty$ and
$\{ f_i \}_1^\infty$, respectively, only takes place if the symmetric
orthogonalization is unique, \cite[Cor.~6]{Tib}. These results will be presented
in section three of the present paper. Complementary, the third author introduced
the notion of a weak symmetric orthogonalization, and he showed its possible
non-uniqueness in all cases by example.

\smallskip 
The concept of frames, i.e.~of sets of generators $\{ f_i \}_{i \in \mathbb N}$
of Hilbert spaces with the property that the inequality
$C \cdot \|x\|^2 \leq \sum_{j \in \mathbb N} |\langle x,f_j \rangle|^2
\leq D \cdot \|x\|^2$ is fulfilled for any $x \in H$ and two constants $C$,
$D>0$, generalizes the notion of a basis for Hilbert spaces. Frames play an
important role in wavelet theory and its applications to signal processing,
image and data compression or analysis, and others (\cite{AF}). Note that in
infinite-dimensional spaces the concepts of frames and bases of Hilbert spaces
do not coincide any more since some bases lack the frame property, see Example
\ref{ex-nonframe}.

The goal of the present paper is to investigate the existence and uniqueness
of symmetric approximations of frames (respectively, symmetric
orthogonalizations of bases of Hilbert spaces) $\{ f_i \}_{i \in \mathbb N}$
of Hilbert subspaces $K \subseteq H$. That means, we look for the existence
and uniqueness of normalized tight frames (resp., orthonormal bases)
$\{ \nu_i \}_{i \in \mathbb N}$ of Hilbert subspaces $L \subseteq H$ such that
the sum $\sum_{j \in \mathbb N} \| \nu_j-f_j \|^2$ is finite and admits the
minimum of all finite sums $\sum_{j \in \mathbb N} \| \mu_j - f_j \|^2$ that
might appear for any other normalized tight frame (resp., orthonormal basis)
$\{ \mu_i \}_{i \in \mathbb N}$ of any other Hilbert subspace of $H$. 
We apply the approach of the third author to symmetric
approximations of frames in Hilbert spaces by normalized tight frames. We can
rely on fundamental work done by David R.~Larson and his collaborators Xingde
Dai, Deguang Han, E.~J.~Ionascu and C.~M.~Pearcy (\cite{DaLa,HaLa,ILP}), by
A.~Aldroubi \cite{Ald}, P.~G.~Casazza \cite{Cas}, O.~Christensen
\cite{Chr1,Chr2,CH} and J.~R.~Holub \cite{Ho1,Ho2}.
As a result we obtain that a symmetric approximation exists and is always unique
if and only if the operator $(P-|F|)$ is Hilbert-Schmidt, where $(I-P)$ denotes
the projection to the kernel of $|F|$. For the detailed explanations we refer to
section two.

\smallskip 
Section one is devoted to an alternative proof of the existence and uniqueness
of the symmetric approximation of finite frames in subspaces of Hilbert spaces
that remodels the proof given by the third author for the case of finite bases.
Some obvious changes appear because of the in general non-trivial kernel of
the operator $F$.

\section{Symmetric approximation of a finite frame in a subspace}

Let $H$ be a Hilbert space and $\{ f_i \}_{i \in \mathbb N} \subset H$ be a
{\it frame} in a separable Hilbert subspace $K \subseteq H$, i.e.~ there are
two constants $C,D > 0$ such that the inequality
\[
  C \cdot \|x\|^2 \leq \sum_{j \in \mathbb N} | \langle x,f_j \rangle |^2 \leq
  D \cdot \|x\|^2
\]
holds for every $x \in K \subseteq H$ and every finite or countable index set.
Without loss of generality we consider only finite and countable frames
in this paper. In case of uncountable frames there are no principal changes.
If $C=D$ then the frame is said to be {\it tight}, if $C=D=1$ then it is
said to be {\it normalized tight}. If for two frames $\{ f_i \}_{i \in
\mathbb N}$ and $\{ g_i \}_{i \in \mathbb N}$ of two Hilbert subspaces $K$ and
$L$ of $H$, respectively, there exists an invertible bounded linear operator
$T: K \to L$ such that $T(f_i)=g_i$ for any index $i$ then these two frames
are said to be {\it weakly similar}. If $K = L$ then the frames are called 
{\it similar}. A {\it Riesz basis} of a Hilbert space $K$ is a 
basis of $K$ that is a frame at the same time. Riesz bases are precisely the
images of orthonormal bases under invertible linear operators,
\cite[Prop.~1.5]{HaLa}. A frame is said to be a {\it near-Riesz basis} if the
deletion of finitely many elements from the frame leads to a set that is a
Riesz basis of the Hilbert space generated by the frame.
P.~G.~Casazza and O.~Christensen constructed in \cite[Prop.~2.4]{CasChr} a
tight frame of a separable Hilbert space that is not a near-Riesz basis.
Other examples have been found by K.~Seip \cite{Seip}. However, any frame of a
finite-dimensional Hilbert space is near-Riesz. For the details of the
separable case we refer to the end of section two.

\begin{definition} \label{symmapprox-fin}
  A normalized tight frame $\{ \nu_i \}_1^n$ in a finite-dimensional Hilbert
  subspace $L \subseteq H$ is said to be a {\it symmetric approximation of a
  finite frame} $\{ f_i \}_1^n$ in a Hilbert subspace $K \subseteq H$ if the
  frames $\{ f_i \}_1^n$ and $\{ \nu_i \}_1^n$ are weakly similar and
  the inequality
  \[
     \sum_{j=1}^n \|\mu_j - f_j\|^2 \geq \sum_{j=1}^n \|\nu_j - f_j\|^2
  \]
  is valid for all normalized tight frames $\{ \mu_i \}_1^n$ in
  Hilbert subspaces of $H$ that are weakly similar to $\{ f_i \}_1^n$.
\end{definition}

If the frame $\{ f_i \}_1^n$ is a system of linearly independent vectors
then the set $\{ \nu_i \}_1^n$ will become a symmetric orthogonalization
of $\{ f_i \}_1^n$, i.e.~an orthonormal system, as will be shown below.
For finite frames $\{ f_i \}_{i=1}^n$ of a Hilbert subspace $K \subseteq H$
we have ${\rm dim}(K) \leq n$ since the reconstruction formula shows the
property of a frame to be a set of generators of $K$. Let $\{ e_i \}_1^n$ be
the standard orthonormal basis of ${\mathbb C}^n$ and define an operator
$F: {\mathbb C}^n \to H$ by the formula
\[
    F\left( \sum_{j=1}^n \alpha_je_j \right) = \sum_{j=1}^n \alpha_jf_j \, .
\]
The operator $F$ has a natural polar decomposition $F=W|F|$, where $W$
is a partial isometry from ${\mathbb C}^n$ into $H$ with initial space
${\rm ker}(F)^\bot \subseteq {\mathbb C}^n$ and range $K \subseteq H$,
cf.~\cite[Prop.~1.10]{HaLa}. In particular, $F$ possesses a closed range and
the set $\{ W(e_i) \}_1^n$ is a normalized tight frame in $K \subseteq H$,
cf.~\cite[Cor.~1.2.(i)]{HaLa}. So the frames $\{ f_i \}_1^n$ and $\{ W(e_i)
\}_1^n$ span the same Hilbert subspace $K \subseteq H$. We want to show that
the normalized tight frame $\{ W(e_i) \}_1^n$ in $K$ is a symmetric
approximation of our initial frame. To proceed we need the following fact
characterizing the Hilbert-Schmidt norm of operators on $H$.

\begin{lemma}    \label{HS-norm}
   Let $T$ be a bounded linear operator on a Hilbert space $H$ and
   let $\{ h_i \}_{i \in {\mathbb J}_1}$ and $\{ k_i \}_{i \in {\mathbb J}_2}$
   be two normalized tight frames of $H$. If one of the sums in the following
   equality is finite then the identity $\sum_{j \in {\mathbb J}_1}\|T(h_j)\|^2
   = \sum_{j \in {\mathbb J}_2}\|T(k_j)\|^2$ holds. If there exists at least
   one normalized tight frame of $H$ for which the sum is finite then $T$ is a
   Hilbert-Schmidt operator and the square root of the sum equals its
   Hilbert-Schmidt norm $\| T \|_{c_2}$.
\end{lemma}

\begin{proof}
   Let $\{ v_i \}_{i \in \mathbb I}$ be an orthonormal basis of the norm-closure
   of the range of $T$ in $H$. Since $T(h_j) \in {\rm ran}(T)$ we have
   \begin{eqnarray*}
      \sum_{j \in {\mathbb J}_1} \|T(h_j)\|^2 & = &
  \sum_{j \in {\mathbb J}_1} \sum_{i \in \mathbb I} |\langle T(h_j),v_i
  \rangle|^2
   = 
  \sum_{j \in {\mathbb J}_1} \sum_{i \in \mathbb I} |\langle h_j,T^*(v_i)
  \rangle|^2 \\
  & = &
  \sum_{j \in {\mathbb J}_1} \sum_{i \in \mathbb I} |\langle T^*(v_i),h_j
  \rangle|^2
   = 
  \sum_{i \in \mathbb I} \sum_{j \in {\mathbb J}_1} |\langle T^*(v_i),h_j
  \rangle|^2 \\
  & = &
  \sum_{i \in \mathbb I}  \| T^*(v_i) \|^2      \,\, .
   \end{eqnarray*}
   The sums are all at most countable (or they are equal to infinity) since a
   Hilbert-Schmidt operator is compact and, hence, its domain and codomain are
   at most separable Hilbert subspaces.
   Since the same calculations are true for the other normalized tight frame
   $\{ k_i \}_{i \in \mathbb I}$ we obtain the desired equality.
\end{proof}

Note that the number of elements in the normalized tight frames, and hence the
dimension of $H$, can be finite or infinite. Moreover, if the frames are finite
then they can contain different numbers of elements. Along the way we have
obtained an interesting new formula for the calculation of the Hilbert-Schmidt
norm of a Hilbert-Schmidt operator. The proof of the next theorem is a
modification of the analogous proof for the symmetric orthogonalization of
finite sets of linearly independent elements $\{ f_i \}_1^n \subset H$ by
the third author, \cite[Th.~1]{Tib}. (An alternative proof for this special case
can be found in \cite{AEG2}.)

\begin{theorem} \label{th-symmapprox-fin}
  Let $\{ \mu_i \}_1^n$ be a normalized tight frame in a Hilbert subspace
  $L \subseteq H$ and let $\{ f_i \}_1^n$ be a frame in a Hilbert subspace
  $K \subseteq H$ such that both these frames are weakly similar. Using the
  notations introduced after Definition \ref{symmapprox-fin} the inequality
  \[
  \sum_{j=1}^n \| \mu_j-f_j \|^2 \geq \sum_{j=1}^n \|W(e_j)-f_j \|^2 =
  \| (I-|F|) \|^2_{c_2} - N
  \]
  holds for $N = {\rm dim}({\rm ker}(|F|))$. Equality appears if and only if
  $\mu_j=W(e_j)$ for every $j =1,2,...,n$. Consequently, the symmetric
  approximation of a frame $\{ f_i \}_1^n$ in a finite-dimensional Hilbert
  space $K \subseteq H$ is a normalized tight frame spanning the same Hilbert
  subspace $L \equiv K$ of $H$ and being similar to $\{ f_i \}_1^n$.
\end{theorem}

\begin{proof}
  Let $\{ \mu_i \}_1^n$ be a normalized tight frame for some Hilbert subspace
  $L \subseteq H$. Define $G: {\mathbb C}^n \to H$ by $G(e_i)=\mu_i$ for
  $1 \leq i \leq n$. Note, that $G$ is a partial isometry by
  \cite[Prop.~1.1]{HaLa}. Since $|F|$ is compact it possesses an orthonormal
  basis of eigenvectors $\{ h_i \}_1^n \subset {\mathbb C}^n$ with eigenvalues
  $\{ \lambda_i \}_1^n$. By Lemma \ref{HS-norm} we can proceed with the
  following calculations:
  \begin{eqnarray*}
    \sum_{j=1}^n \| \mu_j - f_j \|^2 & = &
    \sum_{j=1}^n \| G(e_j)-F(e_j) \|^2 = \sum_{j=1}^n \| (G-W|F|)(e_j) \|^2 \\
    & = &
    \sum_{j=1}^n \| (G-W|F|)(h_j) \|^2
       = \sum_{j=1}^n \| G(h_j)-\lambda_j W(h_j) \|^2  \,\, .
  \end{eqnarray*}
  We point out that the ranges of the (isometric) frame transforms $G^*$ and
  $F^*$ of the frames $\{ \mu_j \}_1^n$ and $\{ f_j \}_1^n$, respectively,
  coincide in ${\mathbb C}^n$ since the frames were supposed to be weakly
  similar, cf.~\cite[Cor.~2.8]{HaLa}. Therefore, the kernels of $G$ and $F$
  also coincide. For eigenvalues $\lambda_i \not= 0$ of $|F|$ we have the
  following lower estimate:
  \begin{eqnarray*}
     \| G(h_i)-\lambda_i W(h_i) \|^2 & = &
        \| G(h_i) \|^2 -2\lambda_i {\rm Re} \langle G(h_i),W(h_i) \rangle +
        \lambda_i^2 \| W(h_i) \|^2 \\
     & = &
        \| G^*G(h_i) \|^2 -2\lambda_i {\rm Re} \langle G(h_i),W(h_i) \rangle +
        \lambda_i^2 \| W(h_i) \|^2 \\ 
     & = & 1-2\lambda_i {\rm Re} \langle G(h_i),W(h_i) \rangle + \lambda_i^2\\
     & \geq & 1-2\lambda_i+\lambda_i^2 = (1-\lambda_i)^2
  \end{eqnarray*}
  for every $i = 1,2, ...,n$ since $h_i \in {\rm ker}(G)^\bot =
  {\rm ker}(W)^\bot$ and
  \[
  {\rm Re} \langle G(h_i),W(h_i) \rangle \leq | \langle G(h_i),W(h_i) \rangle |
   \leq \|G(h_i)\| \|W(h_i)\| = \|G^*G(h_i)\| \|W(h_i)\|=1 \, \, .
  \]
  For eigenvalues $\lambda_i =0$ the term simply equals to zero. As a
  consequence we get the inequality
    \begin{equation}
        \sum_{j=1}^n \| \mu_j-f_j \|^2 \geq \sum_{j=1}^n (1-\lambda_j)^2 - N
    \end{equation}
  valid for $N = {\rm dim}({\rm ker}(|F|))$ and for any normalized tight frame
  $\{ \mu_i \}_1^n$ in a Hilbert subspace $L \subseteq H$ that is weakly similar
  to the frame $\{ f_i \}_1^n$ in $K$. To express the lower estimate in two
  other ways we transform it further using Lemma \ref{HS-norm}:
  \begin{eqnarray*}
     \sum_{j=1}^n (1-\lambda_j)^2 & = &
       \sum_{j=1}^n (1-\lambda_j)^2 \| h_j \|^2 =
       \sum_{j=1}^n \|(1-\lambda_j) h_j \|^2 \\
     & = & \sum_{j=1}^n \|(I-|F|)(h_j) \|^2  = \| (I-|F|) \|_{c_2}^2   \\
     & = & \sum_{j=1}^n \|(I-|F|)(h_j) \|^2  =
       \sum_{j=1}^n \|W(I-|F|)(h_j) \|^2 + N \\
     & = & \sum_{j=1}^n \|W(I-|F|)(e_j) \|^2 + N =
       \sum_{j=1}^n \| W(e_j)-f_j \|^2 + N \, .
  \end{eqnarray*}
  Summing up we get the estimate
  \[
     \sum_{j=1}^n \| \mu_j - f_j \|^2  \geq \sum_{j=1}^n \| W(e_j)-f_j \|^2 =
     \| (I-|F|) \|_{c_2}^2 - N
  \]
  that is valid for every normalized tight frame $\{ \mu_i \}_1^n$ in Hilbert
  subspaces $L$ of $H$ which is weakly similar to the frame $\{ f_i \}_1^n$.

  Finally, we show uniqueness. Suppose, $\{ \mu_i \}_1^n$ is a normalized
  tight frame in a Hilbert subspace $L$ of $H$ that realizes the equality
  $\sum_{j=1}^n \| \mu_j - f_j \|^2 =   \sum_{j=1}^n \| W(e_j)-f_j \|^2$
  and is weakly similar to the frame $\{ f_i \}_1^n$.
  Then ${\rm Re}\langle G(h_i),W(h_i) \rangle =1$ for each $i \in \{1,2,...
  ,n\}$ with $\lambda_i \not=0$. (If $\lambda_i=0$ then $G(h_i)=W(h_i)=0$.)
  Since the inequality
  \[
  1 = {\rm Re} \langle G(h_i), W(h_i) \rangle \leq | \langle G(h_i),W(h_i)
  \rangle | \leq \|G(h_i)\| \|W(h_i)\| =1
  \]
  holds for every eigenvector $h_i$ with non-zero eigenvalue $\lambda_i$
  we conclude $G(h_i)=\alpha_i W(e_i)$ for certain complex numbers
  $\alpha_i \not= 0$ and $i \in \{1,2,...,n\}$ by the Cauchy-Schwarz inequality.
  Replacing $G(h_i)$ we obtain $1=\langle G(h_i),W(h_i) \rangle = \alpha_i
  \langle W(h_i),W(h_i) \rangle = \alpha_i$ for every $i \in \{1,2,...,n\}$.
  Therefore, $G(h_i)= W(h_i)$ for $i \in \{1,2,...,n\}$ forcing $G=W$ since
  $\{ h_i \}_1^n$ has been selected as an orthonormal basis of eigenvectors of
  $|F|$ in $H$. This gives the desired result $\mu_i = W(e_i)$ for every $i \in
  \{1,2,...,n\}$, and the Hilbert subspaces $L$ and $K$ coincide.
\end{proof}

Let us remark that in case $\{ f_i \}_1^n$ is a linearly independent set of
elements in $H$ then $\{ W(e_i) \}_1^n$ has to be linearly independent, too,
since the linear span of both these sets coincides in $H$ and is a linear
subspace of dimension $n$. In this case $|F|$ has no kernel.

\section{Symmetric approximation of frames in separable Hilbert spaces}

Let $H$ be a separable Hilbert space and $\{ f_i \}_{i \in \mathbb N}$ be a
frame in a Hilbert subspace $K \subseteq H$.

\begin{definition} \label{symmapprox-infin}
  A normalized tight frame $\{ \nu_i \}_{i \in \mathbb N}$ in a Hilbert subspace
  $L \subseteq H$ is said to be a {\it symmetric approximation of} $\{ f_i
  \}_{i \in \mathbb N}$ if it is weakly similar to $\{ f_i \}_{i \in \mathbb N}$,
  the inequality
  \[
     \sum_{j=1}^\infty \|\mu_j - f_j\|^2 \geq \sum_{j=1}^\infty \|\nu_j - f_j\|^2
  \]
  is valid for all normalized tight frames $\{ \mu_i \}_{i \in \mathbb N}$ in
  Hilbert subspaces of $H$ that are weakly similar to $\{ f_i \}_{i \in \mathbb N}$,
  and the sum at the right side of this inequality is finite.
\end{definition}

Again we do not assume that both frames $\{ f_i \}_{i \in \mathbb N}$ and
$\{ \nu_i \}_{i \in \mathbb N}$ span the same Hilbert subspace of $H$.
Theorem \ref{th-HS-condition} will show that they do. We note that every set
$\{ \mu_i \}_{i \in \mathbb N}$ that is {\it quadratically close to a frame}
$\{ f_i \}_{i \in \mathbb N}$ of some Hilbert subspace of $H$ (i.e.~for which
$\sum_j \|f_j - g_j \|^2 < \infty$) has to be a frame
of the Hilbert subspace spanned by it, too, cf.~\cite[Th.~3]{CasChr}. However,
it is not automatically weakly similar to the original frame since the
position of possibly existing zero elements in these sequences matters. So we
cannot sharpen our definition on the general level.

Let $\{ e_i \}_{i \in \mathbb N}$ be the standard orthonormal basis
of the separable Hilbert space $l_2$. Given a frame $\{ f_i \}_{i \in \mathbb
N}$ of a Hilbert subspace $K \subseteq H$ we define the operator $F: l_2 \to H$
by the formula
\[
     F \left( \sum_{j=1}^\infty \alpha_je_j \right) = \sum_{j=1}^\infty
     \alpha_j f_j   \, .
\]
The operator $F$ has a natural polar decomposition $F = W |F|$, where $W$ is
a partial isometry from $l_2$ into $H$ with initial space ${\rm ker}(F)^\bot
\subseteq l_2$ and range ${\rm ran}(F)^- \subseteq H$, cf.~\cite[Prop.~1.10]{HaLa}.
In particular, $F$ has closed range and the set $\{ W(e_i) \}_{i \in \mathbb N}$
is a normalized tight frame of a separable Hilbert subspace of $H$,
cf.~\cite[Cor.~1.2(i)]{HaLa}. If we simply try to repeat the steps of our
considerations in the previous section we run into difficulties.

\begin{example}     {\rm
  Let $H=l_2$ with the orthonormal basis $\{ e_i \}_{i \in \mathbb N}$. Set
  $F(e_1)=0$ and $F(e_i)=\alpha_{i-1} e_{i-1}$ for a sequence of complex
  numbers $\{ \alpha_i \}_{i \in \mathbb N}$ and any $i \geq 2$. The set
  $\{ f_i \}_{i \in \mathbb N} = \{ F(e_i) \}_{i \in \mathbb N}$ becomes a
  frame for the Hilbert space $H$ if and only if both
  \[
    \inf_{i \in \mathbb N} \, |\alpha_i|^2 > 0 \,\, ,
    \sup_{i \in \mathbb N} \, |\alpha_i|^2 < +\infty \, .
  \]
  This guarantees the boundedness and surjectivity of $F$. Then
  $|F|(e_1)=0$ and $|F|(e_i) = | \alpha_{i-1} | e_i$ for every $i \geq 2$,
  and consequently $W(e_1)=0$ and $W(e_i)= \alpha_{i-1}/|\alpha_{i-1}| \,
  e_{i-1}$ for any $i \geq 2$. By Lemma \ref{HS-norm} the operator $(I-|F|)$
  is Hilbert-Schmidt if and only if the series
  \[
    {\sum}_{i=1}^\infty \| (I-|F|)(e_i) \|^2 = 1 +
    {\sum}_{i=2}^\infty \| (1-|\alpha_{i-1}|) e_i \|^2 = 1 +
    {\sum}_{i=2}^\infty  (1-|\alpha_{i-1}|)^2
  \]
  converges. As can be easily seen the choice $\alpha_i= {\rm const} \not= 1$
  leads to a infinite Hilbert-Schmidt norm and to a situation in which the
  operator $(I-|F|)$ is definitely not Hilbert-Schmidt. Beside this,
  even if a choice of the coefficients $\{ \alpha_i \}_{i \in \mathbb N}$
  would give rise to a Hilbert-Schmidt operator $(I-|F|)$, the operator $F$
  has still a non-trivial but finite-dimensional kernel by construction. Note
  also, that the frame $\{ f_i \}_{i \in \mathbb N}$ is not a Riesz basis of
  $H$ since it contains the zero vector as its first element.

  \smallskip 
  For other examples the kernel of $F$ can be infinite-dimensional. To see this
  take an orthonormal basis $\{ e_i \}_{i \in \mathbb N}$ of the Hilbert space
  $l_2$  and consider the normalized tight frame $\{ f_i \}_{i \in \mathbb N}$
  for a subspace of $l_2$ defined by $f_{2i}=e_{2i}$, $f_{2i+1}=0$. This frame
  is its own symmetric approximation. However, the kernel of the operator $F$
  is infinite-dimensional since the frame contains countably many zero elements.
  }
\end{example}
  
\begin{theorem} \label{th-HS-condition}
  Let $\{ f_i \}_{i \in \mathbb N}$ be a frame of a Hilbert subspace $K$ of
  the Hilbert space $H$. Denote by $P$ the projection of $H$ onto
  the (norm-closed) range of the operator $F^*F$. Then the operator $(P-|F|)$
  is Hilbert-Schmidt if and only if the sum ${\sum}_{j=1}^\infty \| \mu_j-f_j
  \|^2$ is finite for at least one normalized tight frame $\{ \mu_i \}_{i \in
  \mathbb N}$ of a Hilbert subspace $L$ of $H$ that is weakly similar to
  $\{ f_i \}_{i \in \mathbb N}$. In this situation the estimate
  \[
     {\sum}_{j=1}^\infty \| \mu_j-f_j \|^2 \geq
     {\sum}_{j=1}^\infty \| W(e_j)-f_j \|^2 = \|(P-|F|)\|^2_{c_2}
  \]
  is valid for every normalized tight frame $\{ \mu_i \}_{i \in \mathbb N}$
  of any Hilbert subspace of $H$ that is weakly similar to $\{ f_i \}_{i \in
  \mathbb N}$. (The left sum can be infinite for some choices of $\{ \mu_i
  \}_{i \in \mathbb N}$.)

  Equality appears if and only if $\mu_i =W(e_i)$ for any $i \in \mathbb N$,
  where $W$ is the partial isometry of the polar decomposition $F=W|F|$ and
  $\{ e_i \}_{i \in \mathbb N}$ is the orthonormal basis of $H$ used to define
  the operator $F$. Consequently, the symmetric
  approximation of a frame $\{ f_i \}_{i \in \mathbb N}$ in a Hilbert space
  $K \subseteq H$ is a normalized tight frame spanning the same Hilbert
  subspace $L \equiv K$ of $H$ and being similar to $\{ f_i \}_{i \in \mathbb
  N}$.
\end{theorem}

\begin{proof}
Suppose $(P-|F|)$ is Hilbert-Schmidt. Then it is compact and possesses an
orthonormal set $\{ h_i \}_{i \in \mathbb N}$ of eigenvectors with corresponding
non-zero eigenvalues $\{ \lambda_i \}_{i \in \mathbb N}$. This set forms a
basis of the Hilbert subspace $(I-P)(H)$. Complete this basis of $(I-P)(H)$
to a basis of $H$ by adding a basis $\{ h'_i \}_{i \in \mathbb N}$ of $P(H)$.
Note, that $(P-|F|)(h'_i)=0$ for any index $i$.

Let $\{ \mu_i \}_{i \in \mathbb N}$ be a normalized tight frame of a Hilbert
subspace $L \subseteq H$ that is weakly similar to the frame $\{ f_i \}_{i \in
\mathbb N}$. Define $G:l_2 \to H$ by $G(e_i)=\mu_i$ for $i \in \mathbb N$ and
an orthonormal basis $\{ e_i \}_{i \in \mathbb N}$ of $l_2$. Then $G$ is a
partial isometry with the same kernel as $W$ and $|F|$,
cf.~\cite[Cor.~2.8]{HaLa}. Then by Lemma \ref{HS-norm} we have the following
equality:
  \begin{eqnarray}     \label{eqn-ident}    \nonumber
    \sum_{j=1}^\infty \| \mu_j-f_j \|^2 & = &
      \sum_{j=1}^\infty \| G(e_j)-F(e_j) \|^2 =
      \sum_{j=1}^\infty \| (G-W|F|)(e_j) \|^2 \\ \nonumber
    & = &
      \sum_{j=1}^\infty \| (G-W|F|)(h_j) \|^2 +
      \sum_{k \in \mathbb N} \| (G-W|F|)(h'_k) \|^2 \\
    & = &
      \sum_{j=1}^\infty \| G(h_j) -\lambda_j W(h_j) \|^2 + 0 \,\, .
  \end{eqnarray}
Since $\lambda_i \not= 0$ for any index $i$ we have the following lower
estimates:
  \begin{eqnarray} \label{eqn-ident2} \nonumber
     \| G(h_i)-\lambda_i W(h_i) \|^2 & = &
        \| G(h_i) \|^2 -2\lambda_i {\rm Re} \langle G(h_i),W(h_i) \rangle +
        \lambda_i^2 \| W(h_i) \|^2 \\  \nonumber
     & = &
        \| G^*G(h_i) \|^2 -2\lambda_i {\rm Re} \langle G(h_i),W(h_i) \rangle +
        \lambda_i^2 \| W(h_i) \|^2 \\ \nonumber
     & = & 1-2\lambda_i {\rm Re} \langle G(h_i),W(h_i) \rangle + \lambda_i^2\\
     & \geq & 1-2\lambda_i+\lambda_i^2 = (1-\lambda_i)^2
  \end{eqnarray}
since
  \[
  {\rm Re} \langle G(h_i),W(h_i) \rangle \leq | \langle G(h_i),W(h_i) \rangle |
   \leq \|G(h_i)\| \|W(h_i)\| = \|G^*G(h_i)\| \|W(h_i)\|=1 \, \, .
  \]
Therefore, we get the estimate
  \[
    \sum_{j=1}^\infty \| \mu_j -f_j \|^2 \geq \sum_{j=1}^\infty (1-\lambda_j)^2
    = \sum_{j=1}^\infty \|(P-|F|)(h_j)\|^2 = \| (P-|F|) \|^2_{c_2}
  \]
which is valid for all normalized tight frames $\{ \mu_i \}_{i \in \mathbb N}$
being weakly similar to $\{ f_i \}_{i \in \mathbb N}$.

Obviously, the special choice $\mu_i =W(e_i)$ for $i \in \mathbb N$ gives
the equality
  \[
  \sum_{j=1}^\infty \| W(e_j)-f_j \|^2 = \| (P-|F|) \|^2_{c_2}
  \]
by equality (\ref{eqn-ident}). Uniqueness can be shown in the same way as in
Theorem \ref{th-symmapprox-fin}.

\medskip 
Let $\{ \mu_i \}_{i \in \mathbb N}$ be a normalized tight frame in a Hilbert
subspace $L$ of $H$ that is weakly similar to the frame $\{ f_i \}_{i  \in
\mathbb N}$ and for which the sum $\sum_{j=1}^\infty \| \mu_j - f_j \|^2$ is
finite. Consider the operator $T:l_2 \to H$ defined by $T(e_i)=\mu_i-f_i$
for the fixed orthonormal basis $\{ e_i \}_{i \in \mathbb N}$ of $l_2$.
Then
  \[
     \left\| T \left({\sum}_j \alpha_je_j \right) \right\| \leq
     {\sum}_j | \alpha_j | \cdot \| \mu_j-f_j \| \leq
     \| \{ \alpha_j \}_j \| \sqrt{{\sum}_j \| \mu_j-f_j \|^2}
     < \infty
  \]
for any sequence $\{ \alpha_i \}_i \in l_2$ and, hence, $T$ can be approximated
by finite rank operators in norm on $l_2$. So $T=G-F$ is compact for $G(e_i)=
\mu_i$, $(i \in \mathbb N)$. Counting $F^*F= I-T^*G-G^*T+T^*T$ we see that
$F^*F$ equals the identity operator minus a compact one and, hence, must be
diagonalizable on $l_2$. So $|F|$ is diagonalizable, too.

Let $\{ h_i \}_{i \in \mathbb N}$ be an orthonormal system of eigenvectors
with non-zero
eigenvalues of $|F|$, let $\{ h'_i \}_{i \in \mathbb N}$ be an orthonormal
basis of the kernel of $|F|$. Repeating the calculations (\ref{eqn-ident})
and (\ref{eqn-ident2}) we arrive at
  \[
    \infty > \sum_{j=1}^\infty \| \mu_j-f_j \|^2 \geq
      \| (P-|F|) \|^2_{c_2}
  \]
and the operator $(P-|F|)$ is Hilbert-Schmidt, where $P$ denotes the projection
of $l_2$ onto the kernel of $|F|$.
\end{proof}

\begin{corollary}
  Let $\{ f_i \}_{i \in \mathbb N}$ be a frame of a Hilbert subspace $K$ of
  the Hilbert space $H$. Denote by $P$ the projection of $H$ onto
  the (norm-closed) range of the operator $F^*F$. If the operator $(P-|F|)$
  is not Hilbert-Schmidt then the sum ${\sum}_{j=1}^\infty \| \mu_j-f_j
  \|^2$ is infinite for any normalized tight frame $\{ \mu_i \}_{i \in
  \mathbb N}$ of a Hilbert subspace $L$ of $H$ that is weakly similar to
  $\{ f_i \}_{i \in \mathbb N}$. In other words, there does not exist any
  symmetric approximation of the frame $\{ f_i \}_{i \in \mathbb N}$.
\end{corollary}

\begin{corollary} 
  If a frame $\{ f_i \}_{i \in \mathbb N}$ of an infinite-dimensional Hilbert
  subspace $K$ of a Hilbert space $H$ is a Riesz basis and it has a
  symmetric orthogonalization in $H$ by an orthonormal basis $\{ \nu_i
  \}_{i \in \mathbb N}$ of a Hilbert subspace $L$ of $H$ then the operator
  $(I-|F|)$ is Hilbert-Schmidt and the estimate
  ${\sum}_{j=1}^\infty \| \nu_j-f_j \|^2 < \infty$ holds. Moreover, the more
  general estimate $\sum_{j=1}^\infty \|\mu_j-f_j\|^2 \geq \|(I-|F|)\|^2_{c_2}$
  is valid for any infinite orthonormal basis $\{ \mu_i \}_{i \in \mathbb N}$ of
  any separable Hilbert subspace of $H$ and, in particular, $|F|$ is injective.

  If the operator $(I-|F|)$ is Hilbert-Schmidt then the Riesz basis $\{ f_i
  \}_{i \in\mathbb N}$ admits a unique symmetric approximation by the orthonormal
  basis $\{ W(e_i) \}_{i \in \mathbb N}$ of the same Hilbert subspace of $K
  \subseteq H$.

  Conversely, if in this case the operator $(I-|F|)$ is not Hilbert-Schmidt
  then the sum $\sum_{j=1}^\infty \|\mu_j-f_j\|^2$ is infinite for any
  orthonormal basis $\{ \mu_i \}_{i \in \mathbb N}$ of any separable Hilbert
  subspace of $H$.

\end{corollary}

A reference to Definition \ref{symmapprox-infin} and to Theorem
\ref{th-HS-condition} makes the corollaries obvious.

\smallskip
Recall, that a frame $\{ f_i \}_{i \in \mathbb N}$ is said to be a near-Riesz
basis if there is a finite set $\sigma$ for which the set $\{ f_i \}_{i \in
\mathbb N \setminus \sigma}$ is a Riesz basis of the (separable) Hilbert space
$K \subseteq H$ generated by the initial frame. Near-Riesz bases have a number
of alternative equivalent characterizations:

\smallskip
\begin{tabular}{lp{14.22cm}}
 (i) & the kernel of the derived operator $F:l_2 \to H$ is finite-dimensional
     and has dimension ${\rm card} (\sigma)$, (\cite[Th.~2.4, 3.1]{Ho1}),  \\
 (ii) & the sequence $\{ f_i \}_{i \in \mathbb N}$ is Besselian, i.e.~whenever
     $\sum_{j \in \mathbb N} c_jf_j$ converges in $K$ then $\{ c_i \}_{i \in
     \mathbb N} \in l_2$, (\cite[Th.~2.5]{Ho1}), \\
 (iii) & $\sum_{j \in \mathbb N} c_jf_j$ converges in $K$ if and only if the
       sequence of coefficients $\{ c_i \}_{i \in \mathbb N}$ belongs to $l_2$,
       (\cite[Th.~2.5]{Ho1}), \\  
 (iv) & the sequence $\{ f_i \}_{i \in \mathbb N}$ is unconditional,
      i.e.~whenever $\sum_{j \in \mathbb N} c_jf_j$ converges for a sequence
      $\{ c_i \}_{i \in \mathbb N}$ of numbers then it converges
      unconditionally, (\cite[Th.~3.2]{Ho1} and \cite[Th.~3.1]{CasChr2}).
\end{tabular}

\noindent
In case the frame $\{ f_i \}_{i \in \mathbb N}$ is a near-Riesz basis we can
formulate another corollary the formulation of which already predicts the
result for general bases of infinite-dimensional Hilbert subspaces.

\begin{corollary}
  Let $\{ f_i \}_{i \in \mathbb N}$ be a near-Riesz basis of a Hilbert subspace
  of $H$. Then this frame possesses a symmetric approximation $\{ \nu_i \}_{i
  \in \mathbb N}$ if and only if the operator $(I-|F|)$ is Hilbert-Schmidt.
  In this situation the frame $\{ \nu_i \}_{i \in \mathbb N}$ is also a
  near-Riesz basis and the adjoint operators $F$ and $G$ of the corresponding
  frame transforms of both these frames have kernels of the same dimension
  (i.e.~both these frames have the same excess). In fact, their kernels
  coincide.
\end{corollary}

\begin{proof}
By \cite{Ho1} the operator $F$ has a finite-dimensional kernel if and only if
the frame $\{ f_i \}_{i \in \mathbb N}$ is a near-Riesz basis. However, in case
of a finite-dimensional kernel of $F$ the operator $(I-|F|)$ is Hilbert-Schmidt
if and only if $(P-|F|)$ is. So Theorem \ref{th-HS-condition} implies the
main equivalence of the corollary. The operators $F$ and $W=G$ have the same
kernel since $W^*$ is an isometry and ${\rm ran}F = {\rm ran}W$ as shown above.
Hence, the normalized tight frame $\{ W(e_i) \}_{i \in \mathbb N}$ has to be a
near-Riesz basis, too.
\end{proof}

\begin{corollary}
  If for a frame $\{ f_i \}_{i \in \mathbb N}$ of a Hilbert space $K$
  the operator $(I-|F|)$ is Hilbert-Schmidt then the frame is a near-Riesz
  basis and admits a symmetric approximation.
\end{corollary}

\section{Symmetric orthogonalization of arbitrary linearly independent sets of
         elements}

The considerations below are part of the Ph.~D.~thesis of the third author.
For a different treatment of the infinite situation we refer to \cite{AEG1}.

Let $\{ f_i \}_{i \in \mathbb N}$ be an infinite set of linearly independent
elements of a (separable) Hilbert space $H$ and let $\{ e_i \}_{i \in \mathbb N}$
be an orthonormal basis of $l_2$. Consider the linear operator $F$ defined by
$F(e_i)=f_i$ for $i \in \mathbb N$. Since the set $\{ f_i \}_{i \in \mathbb N}$
may lack the frame property for the Hilbert subspace $K \subset H$ generated
by it the operator can in general be unbounded and can possess a non-trivial
kernel.

\begin{example}    \label{ex-nonframe} {\rm
Let $f_i=e_1+e_i$ for $i \geq 2$ and $f_1=e_1$ for a fixed orthonormal basis
$\{ e_i \}_{i \in \mathbb N}$ of the Hilbert space $l_2$. Then $F$ maps the
element $x_n=\sum_{j=1}^n j^{-1} e_j$ into $F(x_n)=\sum_{j=1}^n j^{-1} e_1 +
\sum_{j=2}^n j^{-1} e_j$. For $n \to \infty$ the sequence $\{ x_n \}_{n \in
\mathbb N}$ converges in norm, whereas the sequence of images $\{ F(x_n) \}_{n
\in \mathbb N}$ diverges to infinity on the multiples of the first basis vector
$e_1$. Therefore, the operator is not everywhere defined on $l_2$ and unbounded.}
\end{example}

\begin{example}  {\rm
To give an example of an operator $F$ with non-trivial kernel consider the
linearly independent set $\{ f_i = e_i - i/(i-1) e_{i-1}, f_1=e_1 \}_{i \in
\mathbb N}$ constructed from an orthonormal basis $\{ e_i \}_{i \in \mathbb N}$
of $l_2$. The derived operator $F: e_i \to f_i$ is bounded since
  \[
    \left\| F \left( {\sum}_{j=1}^\infty \alpha_j e_j \right) \right\|  \leq 
      \left\| {\sum}_{j=1}^\infty \alpha_je_j \right\| +
      \left\| {\sum}_{j=2}^\infty \frac{j}{j-1} \alpha_j e_j \right\| \leq
      3 \cdot \left\| {\sum}_{j=1}^\infty \alpha_j e_j \right\|   \,\, .
  \]
Consider the element $x=\sum_{j=1}^\infty j^{-1} e_j$ and its image
$F(x)$. An easy calculation shows that
  \[
    F \left( {\sum}_{j=1}^n \frac{1}{j} e_j \right)=
    {\sum}_{j=1}^n  \frac{1}{j} f_j = \frac{1}{n} e_n
  \]
forcing $F(x)=0$ for $n \to \infty$. As an immediate consequence we obtain
that the sequence $\{ f_i \}_{i \in \mathbb N}$ is not a frame in the Hilbert
space spanned by it. Note, that the operator $(I-|F|^2)$ is not Hilbert-Schmidt
because Lemma \ref{HS-norm} applied to the basis $\{ e_i \}_{i \in \mathbb N}$
gives an infinite Hilbert-Schmidt norm. Since $(I-|F|^2)=(I-|F|)(I+|F|)$ and
the Hilbert-Schmidt operators form an ideal the operator $(I-|F|)$ also cannot
be Hilbert-Schmidt.  }
\end{example}

To give a reasonable definition of symmetric orthogonalization(s) of infinite
sets of linearly independent vectors $\{ f_i \}_{i \in \mathbb N} \subset K
\subseteq H$ we have to suppose that the derived operator $F$ is at least
bounded. This condition does not depend on the choice of the orthonormal
basis in $l_2$, but only on the set $\{ f_i \}_{i \in \mathbb N} \subset H$.

\begin{definition} \label{symmortho-infin}
  An orthonormal basis $\{ \nu_i \}_{i \in \mathbb N}$ for a Hilbert subspace
  $L \subseteq H$ is said to be a {\it symmetric orthogonalization of} $\{ f_i
  \}_{i \in \mathbb N}$ if the inequality
   \[
   \sum_{j=1}^\infty \|\mu_j - f_j\|^2 \geq \sum_{j=1}^\infty \|\nu_j - f_j\|^2
   \]
  is valid for all orthonormal sets $\{ \mu_i \}_{i \in \mathbb N}$ in
  Hilbert subspaces of $H$ and the sum at the right side of this inequality is
  finite.
\end{definition}

\begin{proposition}  \label{prop-HS-crit}
  Let $\{ f_i \}_{i \in \mathbb N}$ be an infinite set of linearly independent
  elements of a Hilbert space $H$ and let $\{ e_i \}_{i \in \mathbb
  N}$ be an orthonormal basis of $l_2$. Consider the linear operator $F$ defined
  by $F(e_i)=f_i$ for $i \in \mathbb N$. \newline
  If $(I-|F|)$ is Hilbert-Schmidt then  $\sum_{j=1}^\infty \| \mu_j-f_j \|^2
  \geq \| I-|F| \|^2_{c_2}$ for all orthonormal subsets $\{ \mu_i \}_{i \in
  \mathbb N}$ of $H$. \newline
  If $\sum_{j=1}^\infty \| \mu_j-f_j \|^2 < \infty$ for some orthonormal set
  $\{ \mu_i \}_{i \in \mathbb N}$ then the operator $(I-|F|)$ is Hilbert-Schmidt
  and the estimate $\sum_{j=1}^\infty \| \mu_j-f_j \|^2 \geq \| I-|F|
  \|^2_{c_2}$ is valid.
\end{proposition}
  
\begin{proof}
Let $\{ \mu_i \}_{i \in \mathbb N}$ be a countable orthonormal subset of $H$
and $\{ e_i \}_{i \in \mathbb N}$ an orthonormal basis of $l_2$. Define the
operator $G: l_2 \to H$ by $G(e_i)=\mu_i$ for $i \in \mathbb N$. Then $G$ is
an isometry. If $(I-|F|)$ is supposed to be Hilbert-Schmidt then there exists
an orthonormal basis of eigenvectors of the operator $|F|$ in $l_2$ denoted by
$\{ h_i \}_{i \in \mathbb N}$ with corresponding eigenvalues $\{ \lambda_i
\}_{i \in \mathbb N}$. Moreover, since $(I-|F|)$ is supposed to be
Hilbert-Schmidt the operator $(|F|-I)$ is compact and, hence, the operator
$|F|=I+(|F|-I)$ is Fredholm. Therefore, the kernel of $|F|$ has to be
finite-dimensional. Without loss of generality, let $\{ h_1, ...,h_N \}$ be
the eigenvectors corresponding to eigenvalues zero of $|F|$, where $N={\rm
dim}({\rm ker}(|F|)) \geq 0$. By Lemma \ref{HS-norm} we have the following
equality:
  \begin{eqnarray}     \label{eqn-ident3}    \nonumber
    \sum_{j=1}^\infty \| \mu_j-f_j \|^2 & = &
      \sum_{j=1}^\infty \| G(e_j)-F(e_j) \|^2 =
      \sum_{j=1}^\infty \| (G-W|F|)(e_j) \|^2 \\ 
    & = &
      \sum_{j=1}^\infty \| (G-W|F|)(h_j) \|^2 =
      \sum_{j=1}^\infty \| G(h_j) -\lambda_j W(h_j) \|^2  \,\, .
  \end{eqnarray}
Since $G$ is an isometry we have the following lower estimates for every
$i \in \mathbb N$:
  \begin{eqnarray} \label{eqn-ident4} \nonumber
     \| G(h_i)-\lambda_i W(h_i) \|^2 & = &
        \| G(h_i) \|^2 -2\lambda_i {\rm Re} \langle G(h_i),W(h_i) \rangle +
        \lambda_i^2 \| W(h_i) \|^2 \\  \nonumber
     & = & 1-2\lambda_i {\rm Re} \langle G(h_i),W(h_i) \rangle + \lambda_i^2\\
     & \geq & 1-2\lambda_i+\lambda_i^2 = (1-\lambda_i)^2
  \end{eqnarray}
since
  \[
  {\rm Re} \langle G(h_i),W(h_i) \rangle \leq | \langle G(h_i),W(h_i) \rangle |
   \leq \|G(h_i)\| \|W(h_i)\| = 1 \, \, .
  \]
Therefore, we get the estimate
  \[
    \sum_{j=1}^\infty \| \mu_j -f_j \|^2 \geq \sum_{j=1}^\infty (1-\lambda_j)^2
    = \sum_{j=1}^\infty \|(I-|F|)(h_j)\|^2 = \| (I-|F|) \|^2_{c_2}
  \]
which is valid for all orthonormal subsets $\{ \mu_i \}_{i \in \mathbb N}
\subset H$.

\medskip 
Now, suppose the finiteness of the sum $\sum_{j=1}^\infty \| \mu_j -f_j \|^2$
for some orthonormal subset $\{ \mu_i \}_{i \in \mathbb N}$ of $H$.
Consider the operator $T:l_2 \to H$ defined by $T(e_i)=\mu_i-f_i$
for the fixed orthonormal basis $\{ e_i \}_{i \in \mathbb N}$ of $l_2$.
Then
  \[
     \left\| T \left({\sum}_j \alpha_je_j \right) \right\| \leq
     {\sum}_j | \alpha_j | \cdot \| \mu_j-f_j \| \leq
     \| \{ \alpha_j \}_j \| \sqrt{{\sum}_j \| \mu_j-f_j \|^2}
     < \infty
  \]
for any sequence $\{ \alpha_i \}_i \in l_2$ and, hence, $T$ can be approximated
by finite rank operators in norm on $l_2$. So $T=G-F$ is compact for $G(e_i)=
\mu_i$, $(i \in \mathbb N)$. Counting $F^*F= I-T^*G-G^*T+T^*T$ we see that
$F^*F$ equals the identity operator minus a compact one and, hence, must be
diagonalizable on $l_2$. So $|F|$ is diagonalizable, too.

Let $\{ h_i \}_{i \in \mathbb N}$ be an orthonormal system of eigenvectors of
$|F|$ with eigenvalues $\{ \lambda_i \}_{i \in \mathbb N}$. Repeating the
calculations (\ref{eqn-ident3}) and (\ref{eqn-ident4}) we arrive at
  \[
    \infty > \sum_{j=1}^\infty \| \mu_j-f_j \|^2 \geq
      \| (I-|F|) \|^2_{c_2}
  \]
and the operator $(I-|F|)$ is Hilbert-Schmidt. Its Hilbert-Schmidt norm
satisfies the desired inequality.
\end{proof}

\begin{example}  {\rm
We give another example demonstrating that even the condition on
\linebreak[4]
$(I-|F|)$ to be Hilbert-Schmidt does not guarantee that ${\rm ker}(|F|)=
\{ 0 \}$. Let $H=l_2$ with 
  \[
    f_1=e_1 \, , \, f_2= \frac{1}{\sqrt{2}} \, (e_1 + e_2) \, , \,
    f_n= \left( \frac{1}{\sqrt{2}} \right)^{n-1} \negthickspace e_1 -
    \sum_{j=2}^{n-1} \left( \frac{1}{\sqrt{2}} \right)^{n-j+1} \negthickspace
    e_j + \frac{1}{\sqrt{2}} \, e_n \,\,\, {\rm for} \,\, n \geq 3 \, .
  \]
Obviously, the set $\{ f_i \}_{i \in \mathbb N}$ is linearly independent in
$l_2$ by construction. Representing $F$ as an infinite matrix and counting
the entries of the infinite matrix that represents $F^*F$ as scalar products
of column vectors we obtain
\[
(I-|F|^2)(e_1) = - \sum_{j=2}^\infty \left( \frac{1}{\sqrt{2}} \right)^{j-1}
\negthickspace e_j
\,\, , \,
(I-|F|^2)(e_i) = - \left( \frac{1}{\sqrt{2}} \right)^{i-1} \negthickspace
e_1 \,\,\, {\rm for} \,\, i \geq 2 \, .
\]
Applying Lemma \ref{HS-norm} for the orthonormal basis $\{ e_i \}_{i \in
\mathbb N}$ we get a finite value for the Hilbert-Schmidt norm of the
operator $(I-|F|^2)$: $\, \| (I-|F|^2) \|_{c_2}= \sqrt{2}$. Since
\[
\| (I-|F|) \|^2_{c_2} = \sum_{j=1}^\infty (1-\lambda_j)^2 \leq
\sum_{j=1}^\infty (1-\lambda_j)^2 (1+ \lambda_j)^2 = \| (I-|F|^2) \|^2_{c_2}
\]
(where $\{ \lambda_i \}_{i \in \mathbb N}$ are the eigenvalues of $|F|$) we
conclude that the operator $(I-|F|)$ is Hilbert-Schmidt, too.

Consider the element $x = -e_1 + \sum_{j=2}^\infty (\sqrt{2})^{-(j-1)} e_j$
that belongs to $l_2$ since its norm is $\| x \|= \sqrt{2}$. Counting the
value of $|F|^2(x)$ we obtain $x \in {\rm ker}(|F|^2)$. The equality
$\| |F|(x) \|^2 = \langle |F|(x),|F|(x) \rangle = \langle |F|^2(x),x \rangle
= 0$ forces $|F|(x)=0$.

The next theorem demonstrates that $\{ f_i \}_{i \in \mathbb N}$ does not
possess a symmetric orthogonalization. }
\end{example}

\begin{theorem}  \label{th-symmortho}
  Let $\{ f_i \}_{i \in \mathbb N}$ be a linearly independent set of a
  (separable) Hilbert space $H$ such that the derived operator $(I-|F|)$ is
  Hilbert-Schmidt. Then there exists a symmetric orthogonalization $\{ \nu_i
  \}_{i \in \mathbb N}$ of this set in $H$ if and only if ${\rm dim}(({\rm ran}
  F)^\bot) \geq {\rm dim}({\rm ker} F)$. In this case, setting $\nu_i = (V+W)
  (e_i)$ for $i \in \mathbb N$, with an orthonormal basis $\{ e_i \}_{i \in
  \mathbb N}$ of $l_2$ and any partial isometry $V: l_2 \to H$ with initial
  space ${\rm ker}|F|$ and ${\rm ran}V \perp {\rm ran W} = {\rm ran} F$ yields
  a symmetric orthogonalization. Moreover, all symmetric orthogonalizations
  arise in this way.

  The set $\{ f_i \}_{i \in \mathbb N}$ possesses a unique symmetric
  orthogonalization $\{ \nu_i \}_{i \in \mathbb N}$ if and only if
  ${\rm ker} F = \{ 0 \}$, if and only if the sets $\{ \nu_i \}_{i \in \mathbb
  N}$ and $\{ f_i \}_{i \in \mathbb N}$ span the same Hilbert subspace of $H$.
  In this case, $\nu_i = W(e_i)$ for any $i \in \mathbb N$ and $V=0$.
  In other words, the linearly independent set $\{ f_i \}_{i \in \mathbb N}$
  is a frame (and, therefore, a Riesz basis) for the Hilbert subspace it spans
  if and only if ${\rm ker} F = \{ 0 \}$.
\end{theorem}

\begin{proof}
Let us first show that ${\rm ran} F = {\rm ran} W$, i.e.~that ${\rm ran}F$
is actually closed. Since $(I-|F|)$ was supposed to be Hilbert-Schmidt the
operator $|F|$ is Fredholm and his range ${\rm ran}|F|$ is closed by the
definition of Fredholm operators. Let $y_n=F(x_n)$, ($n \in \mathbb N$),
form a Cauchy sequence in ${\rm ran}F$, where $\{ x_n \}_{n \in \mathbb N}
\in l_2$. Then for every $\varepsilon > 0$ there exists a number $N$ such
that $\|y_m-y_n\| < \varepsilon$ for all $m,n \geq N$. Therefore,
  \[
    \varepsilon > \|y_m-y_n\| = \|F(x_m-x_n)\| = \| W |F| (x_m-x_n)\| =
    \| |F| (x_m-x_n)\|
  \]
since $|F|(x_m-x_n) \in ({\rm ker}|F|)^\bot$ belongs to the initial space of
$W$. So the sequence $\{ |F|(x_n) \}_{n \in \mathbb N}$ is a Cauchy sequence
in ${\rm ran}|F|$. However, ${\rm ran}|F|$ is closed and, hence, there exists
a $z \in {\rm ran}|F|$ and a $x \in l_2$ such that $|F|(x)=z=\lim_n |F|(x_n)$.
Then
  \[
    \lim_n y_n = \lim_n F(x_n) = \lim_n W|F|(x_n) = W(\lim_n |F|(x_n)) =
    W|F|(x) = F(x)
  \]
and ${\rm ran}F$ is shown to be closed and to coincide with the range of the
partial isometry $W$.

\smallskip 
Suppose ${\rm dim}(({\rm ran}F)^\bot) \geq {\rm dim}({\rm ker}F)$. If ${\rm
ker}F = \{ 0 \}$ then let $V=0$ on $l_2$. If ${\rm ker}F \not= \{ 0 \}$ then
it has to be finite-dimensional since $|F|$ is Fredholm and ${\rm ker}|F| =
{\rm ker}F$. Suppose, $\{ \xi_i \}_1^n$ is an orthonormal basis of ${\rm
ker}F$. Since ${\rm dim}(({\rm ran}F)^\bot) \geq n$ by assumption we can
select another orthonormal set $\{ \rho_i \}_1^n \subset ({\rm ran}F)^\bot$
and define $V(\xi_i)=\rho_i$ for $i=1,...,n$. Then $V$ is an isometry for
elements of ${\rm ker}|F|$. On $({\rm ker}|F|)^\bot$ we set $V$ to be the
zero operator. Consequently, ${\rm ran}V \perp {\rm ran}W$, $\| (V+W)(x) \|
= \| x \|$ for every $x \in l_2$ and the set $\{ (V+W)(e_i) \}_{i \in \mathbb
N}$ is an orthonormal set for the orthonormal basis $\{ e_i \}_{i \in \mathbb
N}$ of $l_2$. Since $(V+W)|F|(x) = F(x)$ for every $x \in l_2$ we obtain the
equality
  \begin{eqnarray*}
    \sum_{j=1}^\infty \| (V+W)(e_j) -f_j \|^2 & = &
      \sum_{j=1}^\infty \| (V+W)(I-|F|)(e_j) \|^2 \\
    & = & \sum_{j=1}^\infty \| (I-|F|)(e_j) \|^2 = \|(I-|F|)\|^2_{c_2}   \, .
  \end{eqnarray*}
By Proposition \ref{prop-HS-crit} the estimate
  \[
    \sum_{j=1}^\infty \| \mu_j-f_j \|^2 \geq \|(I-|F|)\|^2_{c_2} =
    \sum_{j=1}^\infty \| (V+W)(e_j) -f_j \|^2
  \]
is valid for any orthonormal subset $\{ \mu_i \}_{i \in \mathbb N}$ of $H$.
Therefore, we have found a symmetric orthogonalization of the linearly
independent subset $\{ f_i \}_{i \in \mathbb N} \subset H$ inside $H$.

To show the converse, assume ${\rm dim}(({\rm ran}F)^\bot) < {\rm dim}
({\rm ker}F)$ and the existence of a symmetric orthogonalization $\{ \nu_i
\}_{i \in \mathbb N}$ of the set $\{ f_i \}_{i \in \mathbb N}$ in $H$.
Setting $G(e_i)=\nu_i$ for $i \in \mathbb N$ we obtain an isometry $G:l_2 \to
H$. Consequently,
  \[
    \| G-F \|^2_{c_2} = \sum_{j=1}^\infty  \|(G-F)(e_j) \|^2 =
       \sum_{j=1}^\infty \| \nu_j-f_j \|^2 < \infty
  \]
by the definition of a symmetric orthogonalization. So the operator $(G-F)$
is compact and Hilbert-Schmidt. Since $(I-|F|)$ was assumed to be
Hilbert-Schmidt the kernel of $F$ is finite-dimensional, and so is $({\rm ran}
F)^\bot$ by assumption. Hence, $F$ is Fredholm and $G=F+(G-F)$ is Fredholm, too,
by the compactness of $(G-F)$. Moreover, the indices of $G$ and $F$ coincide
and are greater than zero by assumption. However, ${\rm ind}(G) = {\rm dim}
({\rm ker}G)-{\rm dim}(({\rm ran} G)^\bot) = -{\rm dim}(({\rm ran} G)^\bot)$
since $G$ is an isometry. So we arrive at a contradiction: ${\rm dim}(({\rm
ran} G)^\bot) < 0$. The only possible conclusion is the non-existence of a
symmetric orthogonalization in case ${\rm dim}(({\rm ran}F)^\bot) < {\rm dim}
({\rm ker}F)$. 

\smallskip 
We want to know more about the canonical form of symmetric orthogonalizations
in case at least one exists. First of all, if $\{ \nu_i \}_{i=1}^\infty$ is
another symmetric orthogonalization of the linearly independent set $\{ f_i
\}_{i \in \mathbb N}$ beside the orthonormal set $\{ (V+W)(e_i) \}_{i \in
\mathbb N}$ constructed above we have the relation
  \[
    \| (I-|F|) \|^2_{c_2} = \sum_{j=1}^\infty \| (V+W)(e_j) \|^2 \geq
    \sum_{j=1}^\infty \| \nu_j - f_j \|^2 = \| (I-|F|) \|^2_{c_2}
  \]
by Proposition \ref{prop-HS-crit}. Let $G:l_2 \to H$ again be the isometry
defined by $G(e_i)=\nu_i$ for $i \in \mathbb N$. The goal of the subsequent
considerations is to establish that $(G-W)$ is actually a partial isometry with
initial space ${\rm ker}|F|$ and, hence, $\nu_i = ((G-W)+W)(e_i)$ for $i \in
\mathbb N$, where $(G-W)=V$.

Let $x \in {\rm ker}|F|$. Then $\| (G-W)(x) \|=\| G(x) \| = \|x\|$. So the
operator $(G-W)$ acts on ${\rm ker}|F|$ as an isometry. We want to show that
$G=W$ on $({\rm ker}|F|)^\bot$. Let $\{ h_i \}_{i \in \mathbb N}$ be an
orthonormal set of eigenvectors of $|F|$, let $\{ \lambda_i \}_{i \in
\mathbb N}$ be the corresponding sequence of eigenvalues. In analogy to the
considerations at (\ref{eqn-ident3}) and (\ref{eqn-ident4}) we conclude
  \begin{eqnarray*}
    \| (I-|F|) \|^2_{c_2} = \sum_{j=1}^\infty \| \nu_j - f_j \|^2 & = &
      \sum_{j=1}^\infty \| (G-W|F|)(e_j) \|^2 =
      \sum_{j=1}^\infty \| (G-W|F|)(h_j) \|^2  \\
    & = & \sum_{j=1}^\infty \| G(h_j)-\lambda_j W(h_j) \|^2 \\
    & \geq & \sum_{j=1}^\infty (1-\lambda_j)^2 =
      \sum_{j=1}^\infty \| (I-|F|)(h_j) \|^2 \\
    & = & \| (I-|F|)|^2_{c_2}  \, .
  \end{eqnarray*}
thus, $\sum_{j=1}^\infty \| (G-\lambda_j W)(h_j) \|^2 = \sum_{j=1}^\infty
(1-\lambda_j)^2$ which implies ${\rm Re}(\langle G(h_j),W(h_j) \rangle) =1$
and, consequently, $G(h_j)=W(h_j)$ for any eigenvector $h_j$ with non-zero
eigenvalue $\lambda_j$. Hence, $G=W$ on $({\rm ker}|F|)^\bot$ since the
eigenvectors of $|F|$ with non-zero eigenvalues form a basis of this subspace.

Finally, we verify that ${\rm ran}(G-W) \perp {\rm ran}W$. Let $z_1 \in {\rm
ran}(G-W)$, $z_2 \in {\rm ran}W$. There exist $t_1,t_2 \in l_2$ such that
$(G-W)(t_1)=z_1$, $W(t_2)=z_2$. Furthermore, $t_i=x_i+y_i$ with $x_i \in
{\rm ker}|F|$, $y_i \in ({\rm ker}|F|)^\bot$ and $i=1,2$. Now
  \begin{eqnarray*}
    \langle z_1,z_2 \rangle & = & \langle (G-W)(x_1+y_1),W(x_2+y_2) \rangle
       = \langle G(x_1),W(y_2) \rangle  \\
    & = & \langle G(x_1),G(y_2) \rangle = \langle x_1,y_2 \rangle =0
  \end{eqnarray*}
since $W(x_1)=W(x_2)=0$, $G(y_1)=W(y_1)$, $G(y_2)=W(y_2)$ by the location of
$x_i,y_i$ and the established properties of $W,G$. This shows the orthogonality
of the ranges of $(G-W)$ and $W$.

\smallskip 
To establish uniqueness conditions we first note that ${\rm ker}|F| = \{ 0 \}$
implies $V=0$ and, hence, the uniqueness of the symmetric orthogonalization
$\{ W(e_i) \}_{i \in \mathbb N}$. In this case the sets $\{ \nu_i \}_{i \in
\mathbb N}$ and $\{ f_i \}_{i \in \mathbb N}$ span the same Hilbert subspace
of $H$.                         Suppose now, ${\rm dim}({\rm ker}|F|)=n>0$.
Let $\{ \xi_i \}_1^n$ be an orthonormal basis of ${\rm ker}|F|$. Since the
assumed existence of a symmetric orthogonalization of $\{ f_i \}_{i \in \mathbb
N}$ in $H$ implies ${\rm dim}(({\rm ran}F)^\bot) \geq {\rm dim}({\rm ker}F)$
we can choose an orthonormal set $\{ \rho_i \}_1^n$ in $({\rm ran}F)^\bot$.
Define $V_1:l_2 \to H$ by $V_1(\xi_i)=\rho_i$ for $i=1,...,n$ and $V_1=0$ on
$({\rm ker}|F|)^\bot$, and $V_2:l_2 \to H$ by $V_2=-V_1$.
Then $\nu_i=(V_1+W)(e_i)$ and $\nu_i'=(V_2+W)(e_i)$, ($i \in \mathbb N$), are
symmetric orthogonalizations of $\{ f_i \}_{i \in \mathbb N}$. If we assume
uniqueness of the symmetric orthogonalization then the equalities $\nu_i=\nu_i'$
for $i \in \mathbb N$ lead to the contradiction $\rho_i=-\rho_i =0$ for all
$i \in \mathbb N$. This proves the theorem.
\end{proof}

\noindent
{\bf Acknowledgements:}
The first author is indebted to V.~I.~Paulsen and D.~P.~Blecher for their
invitation to work at the University of Houston in 1998 and for financial
support by an NSF grant.



\end{document}